\documentclass[12pt,longbibliography]{article}
\usepackage[utf8]{inputenc}
\usepackage[T1]{fontenc}
\usepackage[a4paper, margin=2.7cm]{geometry}
\usepackage{amsmath, amsthm, amsfonts, amssymb}
\usepackage{mathrsfs}         \usepackage[colorlinks=true,linkcolor=blue,citecolor=blue,urlcolor=blue,breaklinks]{hyperref}
\usepackage{bbm} 
\usepackage{breakurl}
\usepackage[numbers]{natbib}
\usepackage{url}
\DeclareMathOperator{\ext}{ex}

\newtheorem{thm}{Theorem}
\newtheorem{theorem}{Theorem}

\newtheorem{lem}[thm]{Lemma}

\title{Uniform hypergraphs of girth $6$ and $8$ from generalized polygons}
\author{N.\,G.\ Parvatov}

\AtEndDocument{\bigskip{\footnotesize
  \noindent
  \textsc{N.\,G.\ Parvatov.\\ Russia.}\\
  \textit{E-mail address}: \texttt{\href{mailto:ngparvatov@yandex.ru}{ngparvatov@yandex.ru}}
}}

\begin{document}

\maketitle

\begin{abstract}
Let $\ext_r(N,g)$ be the maximum number of edges in an $r$-uni\-form hypergraph on $N$ vertices with girth at least $g$. We are interested in the asymptotic behavior of this value when $N$ is increasing but parameters $g\in\{6,8\}$ and $r\geq3$ are fixed. It is shown that for some positive constants $c$ and $d$, any integer $r\geq3$ and all sufficiently large integers $N$ the inequalities $\ext_r(N,6)\geq N^{\frac{11}{8}-\frac{c}{\sqrt{\log N}}}$ and $\ext_r(N,8)\geq N^{\frac{11}{9}-\frac{d}{\sqrt{\log N}}}$ hold. 
\end{abstract}

\section{Introduction}
A \emph{hypergraph} is a pair $\mathcal{H}=(V,E)$ of finite sets, where $V$ is the set of vertices, $E$ the set of edges, and each edge from $E$ is a nonempty subset in $V$. (In particular, multiple edges are not allowed.) A vertex $v$ and an edge $e$ are \emph{incident} if $v$ is contained in $e$. A hypergraph is \emph{$r$-uni\-form} if each of its edges is incident to exactly $r$ vertices. $2$-Uni\-form hypergraphs are called \emph{graphs}. The number of edges incident to a vertex is called its \emph{degree}. A hypergraph is \emph{$d$-re\-gu\-lar} if each of its vertices has degree $d$. A hypergraph is called \emph{$k$-par\-ti\-te} if its vertices can be partitioned into $k$ disjoint color classes such that no edge contains more than one vertex in common with any of these classes. A \emph{cyc\-le} of lenght $k\geq2$ in a hypergraph is a pair of sequences of $k$ distinct vertices $v_0,\ldots,v_{k-1}$ and $k$ distinct edges $e_0,\ldots,e_{k-1}$ such that $\{v_i,v_{i+1}\}\subseteq e_i$ for $0\leq i\leq k-1$ providing indices are reduced modulo $k$. The minimum lenght of a cycle in a hypergraph is called its \emph{girth}. Uniform hypergraphs of given girth are used in coding theory, see \cite{XingYuan2018,ShangguanTamo2020,HaymakerTaitTimmons2024}. Applications typically require constructing a hypergraph that has as many edges as possible for a given girth and number of vertices.

Denote by $\ext_r(N,g)$ the maximum number of edges in an $r$-uni\-form $N$-ver\-tex hypergraph of girth not less than $g$. Fixig parameters $g\in\{6,8\}$ and $r\geq3$, we will be interested in the asymptotic behavior of the function $\ext_r(N,g)$ as $N\to\infty$. The following upper and lower bounds are currently known. The upper bound $\ext_r(N,g)=O(N^\frac{3}{2})$ was proved in \citep{BollobasGyory2008} for $r=3$, $g\geq6$ and in \citep{LazebnikVerstraete} for $r\geq3$, $g\geq5$. The upper bound $\ext_r(N,g)=O(N^{1+\frac{1}{\lfloor(g-1)/2\rfloor}})$ was proved in \cite{GyoryLemons2012} for  $r\geq3$, $g\geq4$. It is known from \cite{ShangguanTamo2020} the lower bound $\ext_r(N,g)=\Omega(N^{\frac{g-1}{g-2}}(\log N)^{\frac{1}{g-2}})$ for $r\geq3$, $g\geq6$, which is a logarithmic factor better than the earlier lower bound from \cite{XingYuan2018}. For the case $g=6$ and $r\geq3$, a stronger lower bound $\ext_r(N,6)=\Omega(N^{\frac{4}{3}})$ was recently stated in \cite{HaymakerTaitTimmons2024}. Thus, there are upper bounds $O(N^{\frac{3}{2}})$, $O(N^{\frac{4}{3}})$ and lower bounds $\Omega(N^{\frac{4}{3}})$, $\Omega(N^{\frac{7}{6}}(\log N)^{\frac{1}{6}})$ for $\ext_r(N,6)$, $\ext_r(N,8)$, respectively. Below, lower bounds of the form $\ext_r(N,6)\geq N^{\frac{11}{8}-\frac{c}{\sqrt{\log_2 N}}}$ and $\ext_r(N,8)\geq N^{\frac{11}{9}-\frac{d}{\sqrt{\log_2 N}}}$ are obtained (theorems \ref{th_1} and \ref{th_2}). 

\section{Two observations}
Before describing the constructions, let us make two observations. 

\paragraph*{First observation}
Let $\Gamma=(V_1\sqcup V_2,E)$ be a bipartite graph with color classes $V_1,V_2$ and edge set $E$. For a vertex $v$ in $V_2$, the set $N_v=\{u\in V_1|\{u,v\}\in E\}$ of all vertices in $V_1$ incident to the vertex $v$ is called its \emph{neighborhood}. Suppose no two vertices from $V_2$ have the same nighborhoods. Define the hypergraph $\mathcal{H}_1(\Gamma)=(V_1,E_1)$ where edges are neighbourhoods of vertices from $V_2$. In other words $E_1=\{N_v|N_v\neq\emptyset,v\in V_2\}$. 

We see the hypergraph $\mathcal{H}_1(\Gamma)$ is $r$-uni\-form if in the graph $\Gamma$ all vertices from the color class $V_2$ are of degree $r$. It is $d$-re\-gu\-lar if in the graph $\Gamma$ all vertices from $V_1$ have degree $d$. It is $r$-uni\-form and $d$-re\-gu\-lar simultaneously if in the graph $\Gamma$ vertices from $V_1$ have degree $d$ and vertices from $V_2$ have degree $r$. The bipartite graph $\Gamma$ is called \emph{$(d,r)$-bi\-re\-gu\-lar} in this situation. 

Note that cycles of length $2k$ in $\Gamma$ correspond to cycles of length $k$ in the hypergraph $\mathcal{H}_1(\Gamma)$. This means that the girth of the graph $\Gamma$ is equal to twice the girth of the hypergraph $\mathcal{H}_1(\Gamma)$, which is the first observation. This observation is known from \cite{ErskineTuite2022}. 
\paragraph*{Second observation}
Consider a hypergraph $\mathcal{H}=(V,E)$. Suppose that for each edge $u$ in $E$ there is a hypergraph $\mathcal{H}_u=(V_u,E_u)$ with vertex set $V_u\subseteq u$ and edge set $E_u$. Build the hypergraph $\mathcal{H}^*=(V,E^*)$ where $E^*=\cup_{u\in E}E_u$. In other words the hypergraph $\mathcal{H}^*$ is produced from $\mathcal{H}$ by replacing each of its edges $u$ with edges of the hypergraph $\mathcal{H}_u$. Note that the hypergraph $\mathcal{H}^*$ is $r$-uni\-form if every $\mathcal{H}_u$ is so. 

The hypergraph $\mathcal{H}^*$ has girth at least $g$ in assumption that $\mathcal{H}$ and all $\mathcal{H}_u$ are such. Indead, suppose  
the hypergraph $\mathcal{H}^*$ contains a cycle $C$ of length $k<g$, formed by vertices $v_0,\ldots,v_{k-1}$ and edges $e_0,\ldots,e_{k-1}$. Each edge $e_i$ can be extended to some edge $u_i$ of the hypergraph $\mathcal{H}$. Consider the pair of sequences $v_0,\ldots,v_{k-1}$ and $u_0\ldots,u_{k-1}$ of vertices and edges in $\mathcal{H}$. These sequences do not necessarily form a cycle of length $k$, since edges $u_i$ may repeat, but we can choose a pair of subsequences that form a cycle of length $k'\leq k$. In fact, not all edges $u_i$ are equal, since otherwise $C$ would be a cycle of the hypergraph $\mathcal{H}_u$ for $u=u_1=\ldots=u_k$, contradicting the assumption. If some pair of cyclically successive edges $u_i,u_{i+1}$ are equal, then we can remove one of them, simultaneously removing the common vertex $v_{i+1}$ from the vertex sequence. Thus, one can assume that successive edges are distinct. But then we obtain a cycle of the hypergraph $\mathcal{H}$ with length $k'\leq k$ by choosing a subsequence of distinct edges $u_{i+1},\ldots,u_{i+k'}$ such that $u_{i}=u_{i+k'}$, and the vertices $v_{i+1},\ldots,v_{i+k'}$. The resulting contradiction proves that the hypergraph $\mathcal{H}^*$ has girth at least $g$. 

Thus, $\mathcal{H}^*$ has girth at least $g$ if $\mathcal{H}$ and all $\mathcal{H}_u$ have such girth. This is the second observation. 
\section{The construction}
\paragraph{Basic construction}
Recall that a \emph{generalized $n$-gon} of order $(s,t)$ is a bipartite $(s+1,t+1)$-bi\-re\-gu\-lar graph of girth $2n$ and diameter $n$. It is known, see \cite{Thas1995}, that for each prime power $q$ there exists a generalized $6$-gon (hexagon) of order $(q^3,q)$ with 
\begin{align*}
v(q)&=(1+q)(1+q^4+q^8),\\
b(q)&=(1+q^3)(1+q^4+q^8). 
\end{align*}
vertices in its color classes. By the first observation above, there exists some $(1+q)$-uni\-form and $(1+q^3)$-re\-gu\-lar hypergraph $\Delta_q$ of girth $6$ on $v(q)$ vertices with $b(q)$ edges. 

Further, for each positive integer number $r\leq1+q$, one can produce an $r$-uni\-form  hypergraph of girth at least $6$ on $v(q)$ vertices and $\lfloor\frac{1+q}{r}\rfloor b(q)$ edges. Namely, this hypergraph can be obtained from the hypergraph $\Delta_q$ by replacing each of its edges $e$ with $\lfloor\frac{1+q}{r}\rfloor$ disjoint new $r$-ele\-ment edges included in $e$. If $r|(1+q)$, this hypergraph also is $(1+q^3)$-re\-gu\-lar. The number of its edges is equal asymptotically to $\frac{1}{r}N^{\frac{4}{3}}(1+o(1))$ as the number of vertices $N$ tends to infinity. 

The construction of the hypergraph $\Delta_q$ can be generalized as follows. Given a hypergraph $\Delta_{q'}$ (which is $(1+q')$-uni\-form) for sufficiently large $q'$, its edges can be replaced by isomorphic copies of a smaller hypergraph $\Delta_{q''}$ for some $q''$ such that $1+q'\geq v(q'')$. The resulting hypergraph is $(1+q'')$-uni\-form. Similar replacements can be repeated recursively. As a result, for some $q$, we obtain a $(1+q)$-uni\-form hypergraph of girth at least $6$, which then can be made $r$-uni\-form for any $r\leq1+q$, as described above for $\Delta_q$.

\paragraph{The generalized construction}
Let us describe one of the variants of the generalized construction more precisely. Let $p$ be a prime number, and $m$ a positive integer such that $m\geq2$ and $p^{m-1}\geq5$. Consider the numerical sequence $Q_1,Q_2,\ldots$, defined recursively by the conditions $Q_1=p^m$, and $Q_n=pQ_{n-1}^9$ for $n>1$. We will fix the values $m,p$ until the end of the subsection and will write $Q_{p,m,n}$ instead of $Q_n$ when values $p,m$ are not fixed. It can be verified that
\begin{align*}
Q_n&=p^{9^{n-1}(m+\frac{1}{8})-\frac{1}{8}}
\end{align*}
for all $n\geq1$. 
Note that for all $n>1$ the following two inequalities hold
\begin{align*}
\frac{1+Q_n}{v(Q_{n-1})}>\frac{Q_n}{v(Q_{n-1})}\geq p-1.
\end{align*}
Indead, the first inequality is obwious. For the second inequality we write
\begin{align*}
\frac{Q_n}{v(Q_{n-1})}&=\frac{pQ_{n-1}^9}{v(Q_{n-1})}
=\frac{pv(Q_{n-1})-p(v(Q_{n-1})-Q_{n-1}^9)}{v(Q_{n-1})}=\\
&=p-\frac{p(v(Q_{n-1})-Q_{n-1}^9)}{v(Q_{n-1})}
>p-1
\end{align*}
since
\begin{align*}
Q_{n-1}^9<v(Q_{n-1})&=(1+Q_{n-1})(1+Q_{n-1}^4+Q_{n-1}^8)=\\
&=1+Q_{n-1}^4+Q_{n-1}^8+Q_{n-1}+Q_{n-1}^5+Q_{n-1}^9<5Q_{n-1}^8+Q_{n-1}^9,
\end{align*}
and then
\begin{align*}
\frac{p(v(Q_{n-1})-Q_{n-1}^9)}{v(Q_{n-1})}
<\frac{5pQ_{n-1}^8}{Q_{n-1}^9}=\frac{5p}{Q_{n-1}}\leq\frac{5p}{p^m}=\frac{5}{p^{m-1}}\leq1.
\end{align*}

From the second inequality we also conclude that for all $n\geq1$ we have
\begin{align*}
v(Q_n)\leq\frac{p}{p-1}Q_n^9<pQ_n^9=p\cdot p^{9^n(m+\frac{1}{8})-\frac{9}{8}}<p^{9^n(m+\frac{1}{8})}. 
\end{align*}
Now, we define the sequence of $(1+p^m)$-uni\-form hypergraphs $\Gamma_1,\Gamma_2,\ldots$ (also denoted as $\Gamma_{p,m,1},\Gamma_{p,m,2},\ldots$), where the hypergraph $\Gamma_1=\Delta_{Q_1}=\Delta_{p^m}$ was described earlier, and for each $n>1$ the hypergraph $\Gamma_n$ is obtained from the hypergraph $\Delta_{Q_n}$, which is $(1+Q_n)$-uni\-form on $v(Q_n)$ vertices, inductively by replacing each of its edges with $(p-1)$ disjoint isomorphic copies of the hypergraph $\Gamma_{n-1}$. This is possible due to the first inequality above. Copies of the hypergraph $\Gamma_{n-1}$ contained in different edges of the hypergraph $\Delta_{Q_n}$ have no common edges, since its girth is greater than $2$. Thus, the hypergraph $\Gamma_n$ is $(1+p^m)$-uni\-form on $|V_n|$ vertices and $|E_n|$ edges where 
\begin{align*}
|V_n|&=v(Q_n);\\
|E_n|&=
(p-1)\cdot|E_{n-1}|\cdot b(Q_n)
\geq|E_{n-1}|\cdot b(Q_n)\geq\\
&\geq b(Q_{1})\cdot b(Q_{2})\cdots b(Q_{n})
\geq (Q_1\cdot Q_2\cdots Q_n)^{11}=\\
&=p^{\big(9^0(m+\frac{1}{8})-\frac{1}{8}+9^1(m+\frac{1}{8})-\frac{1}{8}+\cdots+9^{n-1}(m+\frac{1}{8})-\frac{1}{8}\big)\cdot11}=
p^{\big(\frac{9^n-1}{8}(m+\frac{1}{8})-\frac{n}{8}\big)\cdot11}=\\
&=p^{\frac{11}{8}\big((9^n-1)(m+\frac{1}{8})-n\big)}
=p^{\frac{11}{8}\big(9^n(m+\frac{1}{8})-(n+m+\frac{1}{8})\big)}. 
\end{align*}

Finaly, for arbitrary positive integer number $r$ such that $2\leq r\leq1+p^m$, we construct the hypergraph $\Gamma^r_{p,m,n}$ by replacing in the $(1+p^m)$-uni\-form hypergraph $\Gamma_n=\Gamma_{p,m,n}$ each of its edges with $\lfloor\frac{1+p^m}{r}\rfloor$ disjoint new $r$-ele\-ment edges. The hypergraph $\Gamma^r_{p,m,n}$ is $r$-uni\-form. Its girth is at least $6$ by the second observation. It contains no fewer edges than $\Gamma_n$ with the same number of vertices.  Thus, the following lemma holds. 
\begin{lem}\label{lem_1}
Suppose $p$ is a prime number and $m,n,r$ are positive integer numbers such that $p^{m-1}\geq5$ and $2\leq r\leq1+p^m$. There exists an $r$-uni\-form hypergraph $\Gamma^r_{p,m,n}$ of girth at least $6$ with $v(Q_{p,m,n})$ vertices and at least $p^{\frac{11}{8}\big(9^n(m+\frac{1}{8})-(n+m+\frac{1}{8})\big)}$ 
edges.  
$\square$
\end{lem}
It can be seen that the number of edges $|E|$ in the hypergraph $\Gamma=\Gamma^r_{p,m,n}$ is related to the number of vertices $|V|$ as follows:
\begin{align*}
|E|\geq\bigg(p^{9^n(m+\frac{1}{8})}\bigg)^{\frac{11}{8}(1-\varepsilon)}>|V|^{\frac{11}{8}(1-\varepsilon)}
\end{align*}
where $\varepsilon=\varepsilon(m,n)=\frac{n+m+\frac{1}{8}}{9^n(m+\frac{1}{8})}$. It is clear that $0<\varepsilon<1$ for all positive integers $m,n$. 

We further show that by adding isolated vertices (not contained in any edge) in $\Gamma$, we can obtain hypergraphs for which the similar inequalities hold for any sufficiently large number of vertices $|V|$ and sufficiently small $\varepsilon$ depending on $|V|$.  

\paragraph{The lower bound for $\ext_r(N,6)$}Suppose $p$ is any prime and $r\geq2$ an integer number. Choose integer numbers $m^*\geq2,n^*\geq1$ such that 
\begin{align*}
p^{m^*-1}\geq 5, p^{m^*}\geq r-1, 9^{n^*-1}\leq m^*<9^{n^*} 
\end{align*}
and put $N^*=v(Q_{p,m^*,n^*})$. 
For any positive integer $N\geq N^*$, below we describe an $r$-uni\-form hypergraph $\Gamma=(V,E)$ of girth at least $6$ with $N$ vertices and asymptotically $N^{\frac{11}{8}(1+o(1))}$ edges. First, we note the following {\bf property}: for any integer $N\geq N^*$ there exists a pair of integers $m=m_{N},n=n_{N}$, such that 
\begin{align*}
m^*\leq m, n^*\leq n, 9^{n-1}\leq m\leq9^{n+1}\textrm{ and }v(Q_{p,m,n})\leq N<v(Q_{p,m+1,n}). 
\end{align*}
This property can be proved by induction on $N$ using the 
directly verifiable equality 
\begin{align*}
Q_{p,9^{n+1}+1,n}=Q_{p,9^n,n+1}. 
\end{align*}
This equality is true for all positive integers $n$ since, as one can be verified by direct calculation, both its sides are equal to $p^{9^{2n}+\frac{1}{8}(9^n-1)}$. The formulated property is true when $N=N^*$ with the pair $(m,n)=(m^*,n^*)$. If it is true for some 
$N\geq N^*$, then 
\begin{align*}
v(Q_{p,m,n})<N+1\leq v(Q_{p,m+1,n}).
\end{align*} 
We set $m_{N+1}=m$ and $n_{N+1}=n$ if the inequality on the right is strict. Suppose that this is not the case. If $m<9^{n+1}$, then the appropriate values are $m_{N+1}=m+1$ and $n_{N+1}=n$. Otherwise, $m=9^{n+1}$ and 
\begin{align*}
v(Q_{p,m+1,n})=v(Q_{p,9^{n+1}+1,n})=v(Q_{p,9^n,n+1})=N+1.  
\end{align*}
We take $m_{N+1}=9^n$ and $n_{N+1}=n+1$ in this situation. This pair is correct since $m^*<9^{n^*}\leq9^n=m_{N+1}$. Thus, the stated property is true for all $N\geq N^*$.

Now suppose that $N\geq N^*$ and $(m,n)=(m_N,n_N)$. The conditions of lemma \ref{lem_1} are satisfied, so the hypergraph $\Gamma^r_{p,m,n}$ exists. Consider the $N$-ver\-tex hypergraph  $\Gamma$ constructed from it by adding $N-v(Q_{p,m,n})$ new isolated vertices. It is clear that $\Gamma$ is $r$-uni\-form and has girth at least $6$. It has $|V|$ vertices and $|E|$ edges, where 
\begin{align*}
N=|V|&<v(Q_{p,m+1,n})<p^{9^n(m+1+\frac{1}{8})};\\
|E|&\geq
p^{\frac{11}{8}\big(9^n(m+\frac{1}{8})-(n+m+\frac{1}{8})\big)}=
\bigg(p^{9^n(m+1+\frac{1}{8})}\bigg)^{\frac{11}{8}\frac{{9^n(m+\frac{1}{8})-(n+m+\frac{1}{8})}}{9^n(m+1+\frac{1}{8})}}>\\
&>N^{\frac{11}{8}\frac{{9^n(m+\frac{1}{8})-(n+m+\frac{1}{8})}}{9^n(m+1+\frac{1}{8})}}. 
\end{align*}  
From the above inequality for $|V|$ and the inequality $9^{n-1}\leq m$ we have 
\begin{align*}
&\log_p N<9^n\bigg(m+1+\frac{1}{8}\bigg)\leq9m\bigg(m+1+\frac{1}{8}\bigg)<9\bigg(m+1+\frac{1}{8}\bigg)^2,\textrm{ what implies }\\
&m+1+\frac{1}{8}> \frac{1}{3}\sqrt{\log_p N}. 
\end{align*}
Considering the exponent in the inequality for $|E|$ more preciselly, we get
\begin{align*}
\frac{{9^n(m+\frac{1}{8})-(n+m+\frac{1}{8})}}{9^n(m+1+\frac{1}{8})}=
&\frac{{9^n(m+1+\frac{1}{8})-9^n-(n+m+\frac{1}{8})}}{9^n(m+1+\frac{1}{8})}=
\\
=&1-\frac{9^n+n+m+\frac{1}{8}}{9^n(m+1+\frac{1}{8})}\geq1-\frac{33}{\sqrt{\log_p N}}
\end{align*}
since 
\begin{align*}
\frac{9^n+(n+\frac{1}{8})+m}{9^n(m+1+\frac{1}{8})}<
\frac{9^n+9^n+9^{n+1}}{9^n(m+1+\frac{1}{8})}=
\frac{11}{m+1+\frac{1}{8}}\leq\frac{33}{\sqrt{\log_p N}}. 
\end{align*}
Finally, we have 
$|E|\geq N^{\frac{11}{8}\big(1-\frac{33}{\sqrt{\log_p N}}\big)}$. 
Thus the following theorem holds. 

\begin{theorem}\label{th_1}There exists a positive real constant $c>0$ such that for any integer $r\geq 2$ and any sufficiently large integer $N$ the following inequality holds:
\begin{align*}
\ext_r(N,6)\geq N^{\frac{11}{8}-\frac{c}{\sqrt{\log_2 N}}}. 
\end{align*} 
\end{theorem} 
$\square$
\section{Hypergaphs from generalized octagons}
The construction described above can be applied for other bipartite graphs, not just generalized hexagons. Below we apply it to the generalized octagon ($8$-gon) of order $(q^2,q)$ and obtain a lower bound for $\ext_r(N,8)$. We will follow the logic of the previous section, omitting some details. It is known, see \cite{Thas1995}, that for each $q$, which is an odd power of $2$, there exists a Ree-Tits octagon of order $(q^2,q)$ with 
\begin{align*}
v'(q)&=(1+q)(1+q^3+q^6+q^9),\\
b'(q)&=(1+q^2)(1+q^3+q^6+q^9). 
\end{align*}
vertices in its color classes. Following the previous section, we obtain by the first observation an $(1+q)$-uni\-form and $(1+q^2)$-re\-gu\-lar hypergraph $\Delta'_q$ of girth at least $8$ with $v'(q)$ vertices and $b'(q)$ edges.  

We define a sequence of numbers $Q'_1,Q'_2,\ldots$ (also denoted as $Q'_{2,m,1},Q'_{2,m,2},\ldots$) such that $Q'_1=q$ and $Q'_n=2(Q'_{n-1})^{10}$ for $n>1$. Note that all numbers $Q'_n$ are odd powers of $2$. Thus, all hypergraphs $\Delta'_{Q'_n}$ exist. 
Writing $q$ as $q=2^m$ (where $2\nmid m$) we have 
\begin{align*}
Q'_n=2^{10^{n-1}(m+\frac{1}{9})-\frac{1}{9}}\textrm{ for all }n\geq1.
\end{align*} 
Assume additionally that $m\geq5$. Under this assumption the inequalities 
\begin{align*}
\frac{1+Q'_n}{v'(Q'_{n-1})}>\frac{Q'_n}{v'(Q'_{n-1})}\geq1,v'(Q'_n)\leq2(Q'_n)^{10}<2^{10^n(m+\frac{1}{9})}
\end{align*}
hold for all $n>1$ and $n\geq1$ respectively. 

Following the description from the previous section, we construct the $(1+2^m)$-uni\-form hypergraph $\Gamma'_n=\Gamma'_{2,m,n}$ of girth at least $8$ inductively by subsituting isomorphic copies of $\Gamma'_{n-1}$ for edges in $\Delta'_{Q_n}$, assuming that $\Gamma'_1=\Delta'_{Q_1}$. The resulting hypergraph $\Gamma'_n$ has $|V_n|$ vertices and $|E_n|$ edges where now 
\begin{align*}
|V_n|=v'(Q'_n);
|E_n|>2^{\frac{11}{9}(10^n(m+\frac{1}{9})-(n+m+\frac{1}{9}))}. 
\end{align*}
Furthermore, for every positive integer $r\leq1+2^m$, we obtain an $r$-uni\-form hypergraph $\Gamma'^r_{2,m,n}$ from the hypergraph $\Gamma'_n$ (by dividing its edges) on the same vertex set with at least as many edges. Thus, the following lemma holds. 
\begin{lem}\label{lem_2}
For positive integers $m,n,r$, such that $2\nmid m$, $m\geq5$ and $2\leq r\leq1+2^m$, there exists an $r$-uni\-form hypergraph $\Gamma'^r_{2,m,n}$ of girth at least $8$ with $v'(Q'_{2,m,n})$ vertices and at least 
$2^{\frac{11}{9}(10^n(m+\frac{1}{9})-(n+m+\frac{1}{9}))}$ 
edges. 
$\square$
\end{lem}

Next, for given integer $r\geq2$ we choose positive integers $m^*$ and $n^*$, such that 
\begin{align*}
2\nmid m^*, m^*\geq5,r\leq1+2^{m^*}, 10^{n^*-1}-1\leq m^*\leq10^{n^*}-1
\end{align*}
and put $N^*=v'(Q'_{2,m^*,n^*})$. 
For all positive integer number $n$ there is the equality 
\begin{align*}
Q'_{2,10^{n+1}+1,n}=Q'_{2,10^n,n+1}=2^{10^{2n}+\frac{1}{9}(10^n-1)}. 
\end{align*}
We use it to find for arbitrary $N\geq N^*$ a pair of integers $m=m_N,n=n_N$ such that 
\begin{align*}
 2\nmid m,m^*\leq m,n^*\leq n,10^{n-1}-1\leq m\leq10^{n+1}-1\textrm{ and }v'(Q'_{2,m,n})\leq N<v'(Q'_{2,m+2,n}). 
\end{align*}
This is always possible, as we can see by induction. Actually, this is true for $N=N^*$, and providing this is true for some $N\geq N^*$, we have 
\begin{align*}
v'(Q'_{2,m,n})<N+1\leq v'(Q'_{2,m+2,n}). 
\end{align*}
If the second inequality is strict, we set $(m_{N+1},n_{N+1})=(m,n)$. Otherwise, the correct pair is
\begin{align*}
(m_{N+1},n_{N+1})=
\begin{cases}
(m+2,n)&\textrm{ if }m\leq10^{n+1}-3;\\
(10^n-1,n+1)&\textrm{ if }m=10^{n+1}-1.
\end{cases} 
\end{align*}
This is obwious for the first case. For the second case we have inequalities 
\begin{align*}
v'(Q'_{2,10^n-1,n+1})<v'(Q'_{2,10^n,n+1})=v'(Q'_{2,10^{n+1}+1,n})=N+1<v'(Q'_{2,10^n+1,n+1}). 
\end{align*}
and $m^*\leq10^{n^*}-1\leq10^{n}-1= m_{N+1}$. 
Thus, a correct pair exists for each $N\geq N^*$. 

Suppose $N\geq N^*$ and $(m,n)=(m_{N},n_{N})$. By Lemma \ref{lem_2}, the hypergraph $\Gamma'^r_{2,m,n}$ exists. Adding to it isolated vertices, we obtain a hypergraph $\Gamma'$ on $|V|$ vertices and $|E|$ edges, where
\begin{align*}
N=|V|<2^{10^n(m+2+\frac{1}{9})}, 
|E|>N^{\frac{11}{9}\frac{10^n(m+\frac{1}{9})-(n+m+\frac{1}{9})}{10^n(m+2+\frac{1}{9})}}
>N^{\frac{11}{9}\big(1-13\sqrt{\frac{10}{\log_2 N}}\big)}, 
\end{align*}  
which can be proven by calculations similar to those in the previous section. 

Thus the following theorem is true. 

\begin{theorem}\label{th_2}There exists a positive real constant $d>0$ such that for any integer $r\geq 2$ and any sufficiently large integer $N$ the following inequality holds:
\begin{align*}
\ext_r(N,8)\geq N^{\frac{11}{9}-\frac{d}{\sqrt{\log_2 N}}}. 
\end{align*} 
\end{theorem} 
$\square$

\bibliography{refs} 

\end{document}